\documentclass[a4paper]{amsart}

\title[Combinatorial aspects of Connes's embedding conjecture]{Combinatorial aspects of Connes's embedding conjecture\\ and asymptotic
distribution\\ of traces of products of unitaries}
\thanks{Research partially supported by NSF GRANT no. DMS 0200741}
\author{Florin R\u{a}dulescu}
\address{Department of Mathematics \\
University of Roma ``Tor Vergata''\\ Via della Ricerca Scientifica, 00133
Roma, Italy\\
on leave from University of Iowa\\ Iowa City, IA 52242, USA}

\usepackage{amsmath,amsthm}
\usepackage{amsmath}
\usepackage{amsfonts}
\usepackage{amssymb}
\usepackage{amsthm}
\usepackage{calc}
\usepackage{times}

\def\tilde{\widetilde}
\def\tr{\mathop{\rm tr}}
\def\Tr{\mathop{\rm Tr}}
\def\id{{\rm Id}}

\def\leq{\leqslant}
\def\geq{\geqslant}
\def\d{\,{\rm d}}
\def\eqno{\leqno}

\def\card{\mathop{\rm card}}
\def\epsilon{\varepsilon}

\def\cU{{\cal U}}
\def\C{{\mathbb C}}
\def\N{{\mathbb N}}
\def\diez{\sharp}

\theoremstyle{plain} 
\newtheorem{thm}{Theorem}[section]

\newtheorem{prop}[thm]{Proposition}

\theoremstyle{definition}
\newtheorem{defn}[thm]{Definition}

\newtheorem{rem}[thm]{Remark}


\def\cal{\mathcal}
\begin{document}


\begin{abstract}
In this  paper we study the asymptotic distribution of the moments of
(non-normalized) traces $\Tr (w_1), \Tr(w_2), \ldots,
\Tr(w_r)$, where $ w_1, w_2, \ldots, w_r$ are reduced
words in unitaries in the group $\cU(N)$. We prove that as $N\to \infty$ these
variables are distributed as normal gaussian variables $\sqrt {j_1} Z_1,
\ldots, \sqrt{Z_r}$, where $j_1, \ldots, j_r$ are the number of cyclic
rotations of the words $w_1, \ldots, w_s$ leaving them invariant. This extends
a previous result by Diaconis (\cite{Diac}), where this it was proved, that
$\Tr(U), \Tr(U^2), \ldots,$ $\Tr(U^p)$ are asymptotically distributed as $Z_1,
\sqrt 2 Z_2, \ldots, \sqrt p Z_p$. 

We establish a combinatorial  formula for $\int |\Tr (w_1)|^2\cdots|
\Tr(w_p)|^2$. In our computation we reprove some results from \cite{BC}.
\end{abstract}

\maketitle

\thispagestyle{empty}

\section{Introduction}

Connes's embedding conjecture (\cite{Co}) for the case of discrete groups states that
every discrete group $\Gamma$ can be asymptotically embedded in the algebra
of $N$ by $N$ matrices, when $N$ tends to infinity.
As  observed in \cite{Ra} (see also \cite{Oz} and \cite {Ki}) it amounts to
prove that for every finite subset $F$ of $\Gamma$, for every $\varepsilon >0$,
there exist $N$ and unitaries $\{ a_f \mid f\in F\}$ in $\cU (N)$ such that
$\| a_{f_1} a_{f_2}- a_{f_1f_2}\|_{\rm HS}\leq \varepsilon \|\id \|_{\rm HS} 
$ and  $f_1f_2\in F$ for all $f_1, f_2$ in $F$. Here by
$\|\cdot\|_{\rm HS}$ we denote the Hilbert-Schmidt norm
\begin{equation*}
\|A\|_{\rm HS} = \Tr (A^* A)^{1/2}, \quad A \in M_N(\C),\end{equation*} 
$\Tr$ being the (non-normalized) trace on $M_N(\C)$. If $\Gamma$ is a group
with presentation $\langle F_\infty\mid R\rangle$, where $R$ are the relators,
it can be proved (see \cite{Ra}) that the Connes's embedding conjecture is
equivalent 
to show that for any $\varepsilon>0$, $w_1, w_2, \ldots, w_s\in R$,
and for any $w_0\not \in R$, assuming that $w_0, w_1, \ldots, w_s$ are the
words on the letters $a_1,\ldots, a_M$, there exist $N$ and unitaries
$U_1, U_2, \ldots, U_P$ in $\cU(N)$ such that if $W_0, \ldots,W_s$ are the
corresponding words obtained by substituting $a_1, \ldots, a_M$ with $U_1,
\ldots, U_p$ we have (with $\Tr =\frac 1 n \Tr$)
\begin{equation}
|\tr (W_0)| <\varepsilon, \quad |\tr (W_1)|>1-\varepsilon, \quad \ldots, \quad
 |\tr (W_s)|> 1-\varepsilon  .\label{stea}
\end{equation}
Consequently, a natural object to study is the following: Let $F_M$ be the
 free group with $M$ generators $a_1, a_2, \ldots, a_M$. Let $w_0, w_1,
 \ldots, w_s$ be the reduced words in $F_M$ and let $f_{w_0}, \ldots,
 f_{w_s}$ be the functions on $(\cU(N))^M$ obtained by evaluating the traces
 $\Tr (W_0), \ldots, \Tr (W_s)$ of the words $W_0, \ldots, W_s$ at an $M$-uple
 $(U_1, \ldots, U_M)$. Then one has to determine the joint moments of these
 functions, i.e.\ the quantities (for all $\alpha _0, \ldots, \alpha_s$ in~ $\N$)
\begin{equation*} \int_{(\cU(N))^M} |f_{w_0}|^{\alpha_1} \cdots
     |f_{w_s}|^{\alpha_s}\d U_1 \cdots \d U_M,\end{equation*}
with respect to the Haar measure.

In particular, after normalizing with the factor $\frac 1N$, if one determines
the measure for these moments, then one could solve the
inequality~\eqref{stea}.

\section{Computation of $\int_{\cU(N)} u_{i_1 j_1} \cdots u_{i_p j_p}
  u^*_{r_1s_1} \cdots u^*_{r_p s_p}\d U$}

The following computation was first performed by D. Weingarten, F. Xu,
and B.~Collins. At the time of writting the paper we were not aware of the
previous literature, so we include our own proof for this computation. 

Let $S_n$ be the group of $n$-permutations and let $\C[S_n]$ be the group
algebra. As in \cite{CS}, we denote by $W^N_\sigma$ the coefficient of $\sigma \in
S_n$ in the inverse of the element $\Phi^N = \sum\limits_{\sigma \in S_n}
N^{\diez\sigma}\sigma \in \C[S_n]$
($\diez \sigma$ is the number of cycles in $\sigma$; 
the element $\Phi ^N$ is invertible as it will be proven bellow for
$N>n$). Thus we take
\begin{equation*}
(\Phi^N)^{-1}= \sum _{\sigma\in S_n} W_\sigma^N \cdot \sigma.\end{equation*}
Note that $\Phi^N$ is a central element and hence so is $\sum\limits_{\sigma\in
  S_n}W^N_\sigma \cdot \sigma.$ With these notations we have:
\begin{thm} 
For $N, n$ in $\N$, $N>n$ and $\d U$ the Haar measure on $\cU(N)$, let $i_1,
\ldots, i_n$, $j_1, \ldots, j_n, r_1, \ldots, j_n, s_1, \ldots, s_n$ be indices
from $1$ to $N$. Denote the entries  of a unitary by $u_{ij}$ and the
entries of its  adjoint by $u^*_{ij}= \overline{u_{ij}}$.
Then 
\begin{equation*}
\int_{\cU(N)} u_{i_1 r_1} \cdots u_{i_n r_n} u^*_{s_1 j_1}\cdots
u^*_{s_nj_n}\d U= \sum W^N_{\sigma \theta^{-1}}\end{equation*}
with the sum in the right hand side running over all $\sigma, \theta $ in
$S_n$ such that $j_a = i_{\sigma(a)}$, $a=1,2, \ldots, n$ and $s_b=
r_{\theta(b)}$, $B=1,2, \ldots,n$.  
\end{thm}
\begin{proof}
Let $L^2(M_N(\C)^n, \mu^n_M)$ be the Hilbert space obtained by endowing 
$M_N(\C)^n$ with the measure $C {\rm e}^{-\Tr (A^*_1 A_1)- \cdots-
  \Tr(A^*_nA_n)}$, $ (A_1, \ldots, A_n)\in M_N(\C)^n$, where $C$ is a constant,
so that the entries functions $(A_1, A_2, \ldots, A_n) \mapsto
a_{i_1j_1}^{(1)}\cdots a_{i_nj_n}^{(n)}$, have norm $1$. Here 
$a_{ij}^{(t)}$ are
the $ij$-entries of the matrix $A^{(t)}$ on the $t$-th component of the
product $(M_N(\C))^n$. 

Denote, for $\sigma $ in $S_n$, by $\chi_\sigma$ the function
\begin{equation*}
\chi_{\sigma} = \sum _{i_1, \ldots, i_n=1}^N a^{(1)}_{i_1 \sigma (i_1)}
a^{(2)}_{i_2 \sigma(i_2)} \cdots a^{(n)}_{i_n\sigma(i_n)}.
\eqno(?)
\end{equation*}

Then, from the theory of symmetric functions 
(\cite{McD}, \cite{Procesi}), the functions $\chi _\sigma$
generate the subspace functions on $(M_N(\C))^n$ 
that are invariant to the diagonal
action of $\cU(n)$
 on $(M_N(\C))^n$: $(A_1, \ldots, A_n)\mapsto (UA_1U^*, \ldots, UA_nU^*)$,
 $U\in \cU(N)$. Moreover, for $n<N$ the functions $\{\chi_\sigma\mid \sigma \in
 S_n\}$ are independent (\cite{Procesi}) 
and the scalar product $\langle \chi_\sigma,
 \chi_\mu\rangle$ depends only on $\sigma^{-1}\mu$ and it is equal to
 $N^{\diez (\sigma^{-1}\mu)}$.

Consequently, $\langle \chi_\sigma , \chi_\mu\rangle_{\sigma, \mu\in S_n}$
represents the matrix of the convolution with $\Phi_N$ on $L^2(S_n)$. 
Consequently, the
inverse of $\Phi_N$ (which exists since the functions are independent) is the matrix
$(W_{\sigma^{-1}\mu})_{\sigma, \mu \in S_n}$. Let $P$ be the projection from
$L^2((M_N(\C))^n), \mu)$ onto the space of $\cU(N)$ invariant functions.
Then on one hand, since $\mu$ is an invariant measure, it follows that $P$ is
the average over $\cU(N)$ by integration. Hence 
\begin{equation}\label{unu}
P(a^{(1)}_{i_1j_1}\cdots a_{i_nj_n}^{(n)}) = \int_{\cU(N)}
(ua^{(1)}u^*)_{i_1j_1}\cdots (ua^{(n)}u^*)_{i_nj_n} \d U.
\end{equation}
On the other hand, assume
$P(a^{(1)}_{i_1j_1}\cdots a_{i_nj_n}^{(n)})\!=\!\! \sum\limits_{\sigma\in S_n} c_\sigma
\chi_\sigma$, where $c_\sigma$ depends on $i_1, \ldots, i_n$, $j_1, \ldots,
j_n$.
Then, for all $\mu\in S_n$,
\begin{equation*}
\Big\langle a_{i_1, j_1}^{(1)} \cdots a_{i_nj_n}^{(n)}-\sum_{\sigma} c_\sigma
\chi_\sigma , \chi_\mu\Big\rangle=0, \end{equation*}
and hence
\begin{equation}\label{doi}
\langle a^{(1)}_{i_1j_1}\cdots a_{i_nj_n}^{(n)}, \chi _\mu\rangle =\sum
_\sigma c_\sigma \langle \chi_\sigma ,\chi_\mu \rangle, \quad \forall \mu \in
S_n\end{equation}
But $\langle a^{(1)}_{i_1j_1}\cdots a_{i_nj_n}^{(n)}, \chi_\mu\rangle_\mu$  is
  the vector (indexed) by $\mu\in S_n$) with the property that the $\mu$-th component
  is equal to $1$ if and only if $j_a = i_{\mu(a)}$, $a= 1, 2, \ldots, n$.

Let $R_{\sigma, \mu}$ be the inverse of the matrix $(\langle \sigma,
\mu\rangle)_{\sigma, \mu\in S_n}$. We have noted before that $R_{\sigma,
  \mu}=W^N_{\sigma^{-1} \mu}$. From~\eqref{doi}, by inversion, we deduce that
$c_\sigma =\sum_\mu R_{\sigma, \mu}$, where the sum runs over all $\mu$ such
that $j_a=i_{\mu(a)}$, $a= 1, 2, \ldots  n$.
 Thus
\begin{equation} \label{trei}
P(a^{(1)}_{i_1j_1}\cdots a_{i_nj_n}^{(n)})\!=\! \sum_{{\scriptstyle \sigma \in
    S_n}\atop{\scriptstyle \mu \in S'}} R_{\sigma, \mu}\chi_\sigma, \quad
    S'= \{\mu\in S_n\mid j_a=i_{\mu (a)}, a=1, 2, \ldots, n\}. 
\end{equation} 
From \eqref{unu} we obtain that 
\begin{align*}
P(a^{(1)}_{i_1j_1}\cdots a_{i_nj_n}^{(n)})&= \sum_{{\scriptstyle r_1, \ldots,
    r_n=1}\atop{\scriptstyle s_1, \ldots, s_n=1}} \int _{\cU(N)} u_{i_1
    r_1}\cdots u_{i_nr_n} a^{(1)}_{r_1s_1}\cdots a_{r_ns_n}^{(n)}
    u^*_{s_1j_1}\cdots u^*_{s_nj_n}\d U\\&= \sum_ {{\scriptstyle r_1, \ldots,
    r_n=1}\atop{\scriptstyle s_1, \ldots, s_n=1}} 
a^{(1)}_{r_1s_1} \cdots a^{(n)}_{r_ns_n}\int _{\cU(N)} u_{i_1r_1}\cdots
    u_{i_nr_n}u^*_{s_1j_1}\cdots u^*_{s_nj_n}\d U.\end{align*}
Identifying the coefficients from the last formula with~\eqref{trei} we obtain
    our statement.
\end{proof}

\section{Formula for $\iint _{\cU(N)^2}\Tr (W_1) \cdots \Tr(W_s)\d U \d V$}

In this section we deduce a formula for the integral of 
traces of words only in case of
$\cU(N)^2$ (instead of $\cU(N)^M$) for simplicity. 
A similar formula was derived in \cite{BC}. Since the shape of the
combinatorial aspect of the  formula is important for the computation of the asymptotics, we derive our formula directly from
the preceding section.
 
Let $w_1, w_2, \ldots, w_n$ be reduced words in $F_2=\langle a, b \rangle$, and   let $W_1, W_2, \ldots, W_n$ be the corresponding words viewed as functions in
the variables $(U, V)\in \cU(N)^2$ obtained by substituting $(U, V)$ for $(a,
b)$.  We describe $\Tr(w_1) \cdots \Tr(w_n)$ in terms of a permutation $\gamma$
and write $\Tr (w_1)\cdots \Tr(w_n)= \Tr _\gamma(w_1\cdots w_n)$, where
$\gamma$ is described as follows.

Let $n$ be the total number of occurrences of the symbol $u$ in $W_1, W_2,
\ldots, W_n$. Let $m$ be the total number of occurrences of $v$. For the
integral $\iint_{\cU(N)} \Tr (W_1)\cdots \Tr(W_p) \d U \d V$ 
to be non-zero it is 
necessary (\cite{Petz}) that $n$ equals the number of $U^*$ and that $m$
equals the numbers of $V^*$. We introduce a set of symbols (indexed by the
letters $u, v, u^*, v^*$ respectively)
\begin{equation*}
X=\{1_u, \ldots, n_u, 1_{u^*}, \ldots, n_{u^*}, 1_v, \ldots, m_v, 1_{v^*},
\ldots, n_{v^*}\}.\end{equation*} 
Thus $X$ is a set with $2(n+m)$ elements.

\begin{defn} Given $w_1, \ldots, w_p$ and $X$ as above we define a
  permutation $\gamma$ of $X$ by means of the formula
\begin{equation*}
\begin{split}
\Tr\nolimits _\gamma(w_1 \cdots w_p)&= \Tr (w_1) \cdots \Tr (w_p)\\&= 
\sum_{{{{\scriptstyle a_{1_u}, \ldots, a_{n_u}=1}\atop {\scriptstyle a_{1_{u^*}},
    \ldots, a_{n_{u^*}}=1}}\atop {\scriptstyle a_{1_v}, \ldots,
  a_{m_{v}}=1}}\atop
{\scriptstyle a_{1_{v^*}}, \ldots, a_{m_{v^*}}=1}}
^N u_{a_{1_u} a_{\gamma(1_u)}}\cdots u_{a_{n_u} a_{\gamma(n_u)}}
u^*_{1_{u^*}\gamma(1_{u^*})}\cdots u^*_{n_{u^*} \gamma(n_{u^*})}
\\&\qquad\qquad\qquad\qquad \cdots v_{1_v \gamma(1_v)} \cdots v^*_{m_{v^*}\gamma(m_{v^*})}
\d u \d v.
\end{split}
\end{equation*}
We denote the term on right hand side by $\Phi_\gamma (U,V)$, where $U,V\in \cU(N)$.
\end{defn}
Note that since the words are reduced $\gamma$ has no fixed points.

With these notations we have:

\begin{thm} The integral of $\Tr (w_1)\cdots \Tr(w_p)$ over $U(N)^2$ is 
\begin{equation*}
\iint_{\cU(N)^2} \Tr\nolimits_\gamma (U, V) \d U \d V=
\sum_{{\scriptstyle \sigma_u, \theta_u\in S(1_u, \ldots, n_u)}\atop{\scriptstyle \sigma
    _v ,\theta_v\in S(1_v, \ldots, n_v)}} W^N_{\sigma_u\circ \theta^{-1}_u}
W^N{\sigma_v \circ \theta^{-1}_v} N^{\diez R(\gamma, \sigma_u, \ldots, \theta _v},\end{equation*}
where $R(\gamma, \sigma_u,\theta_u, \sigma_v, \theta_v)$ is the equivalence
relation on $X$ generated by
\begin{alignat*}{3}
t_{u^*}= \gamma(\sigma_u(t_u)),&\quad& \gamma (t_{u^*}) = \theta _u
(t_u),&\quad& t=1, 2, \ldots, n,\\
s_{v^*}= \gamma(\sigma_v(s_v)),&\quad& \gamma (s_{v^*}) = \theta _v
(t_v),&\quad& s=1, 2, \ldots, m.\!\!\end{alignat*} 
Here $\diez R(\gamma, \sigma_u, \theta_u, \sigma_v, \theta_v)$ is the number
of classes in the equivalence relation.
\end{thm}

\begin{proof}
Indeed, we have that 
\begin{align*} &\iint_{\cU(N)^2} \Tr\nolimits_\gamma (U, V) \d U \d V\\&=
\bigg(\sum_{{\scriptstyle a_{1_u}, \ldots, a_{n_u}=1}\atop{\scriptstyle
    a_{1_{u^*}}, \ldots, a_{n_{u^*}}=1}} \int _{\cU(N)} u_{a_{1_u}
    a_{\gamma(1_u)}}\cdots  u_{a_{n_u}
    a_{\gamma(n_u)}}           u^*_{a_{1_{u^*}}a_{\gamma(1_{u^*})}}\cdots 
u^*_{a_{n_{u^*}}
    a_{\gamma(n_{u^*})}} \d U\bigg)\\
&\qquad \cdot\bigg (\sum_{{\scriptstyle a_{1_v}, \ldots, a_{n_v}=1}\atop{\scriptstyle
    a_{1_{v^*}}, \ldots, a_{n_{v^*}}=1}} \int _{\cU(N)} 
v_{a_{1_v}
    a_{\gamma(1_v)}} \cdots v_{a_{n_v}
    a_{\gamma(n_v)}}    v^*_{a_{1_{v^*}}
    a_{\gamma(1_{v^*})}}\cdots v^*_{a_{n_{v^*}}
    a_{\gamma(n_{v^*})}}\d V\bigg)
\end{align*}
We now apply the formula from the preceding section, interchange the summation
formula with the summation after $\sigma_u, \theta_u, \sigma_v, \theta_v$,
where $\sigma_u, \theta_u$ are the permutations that appear in the integrals
for the $u$'s and $\sigma_v, \theta_v$ are the permutations that appear in the
summations for the $\theta$'s.
\end{proof}
\begin{rem} Since the words are allays reduced, the equivalence relation\break
  $R(\gamma, \sigma _u, \theta_u, \sigma_v, \theta_v)$ has no singleton
  classes and hence $\diez R (\gamma, \sigma_u, \theta_u,
  \sigma_v,\theta_v)\leq n+m$.
\end{rem}

\section{The asymptotics for $\iint (\Tr
  (W_1))^{\alpha_1}(\overline{\Tr(W_1)})^{\beta_1}\cdots (\Tr(W_p))^{\alpha_p}(\overline{\Tr
  (W_P)})^{\beta_p} \d U \d V$}

In this section we show that for all words $w_1, \ldots, w_p$ in $F_2$ 
by taking
$W_1, \ldots, W_p$ to be the corresponding functions on $\cU(N)^2$ we have that 
\begin{equation*}\iint_{\cU(N)} (\Tr
  (W_1))^{\alpha_1}(\overline{\Tr(W_1)})^{\beta_1}\cdots (\Tr(W_p))^{\alpha_p}
(\overline{\Tr
  (W_P)})^{\beta_p} \d U \d V\end{equation*}
 is ${\rm O}\big(\frac 1N\big)$ unless
  $\alpha_1=\beta_1, \ldots, \alpha_p=\beta_p$ in which case the integral is
  $\alpha_1! \cdots\alpha_p!$\break $\cdot(j(w_1))^{\alpha_1} \cdots (j(w_p))^{\alpha^p}$,
  where $j_i=j(w_i)$ is the numbers of cyclic rotations of the word $w_i$
  which leave $w_i$ it invariant. 
As in \cite{Diac} this means that the asymptotic
  distribution of the variables $\Tr (W_1), \ldots, \Tr(W_p)$ as $M\to\infty$
  is that of $\sqrt {j_1} Z_1 , \ldots, \sqrt{j_p} Z_p$, where $Z_1, \ldots,
  Z_p$ are independent gaussian variables.
\begin{thm} Let $w_1, \ldots, w_p$ be the words on $F_2$ and $W_1, \ldots,
  W_p$ be the corresponding functions on $\cU(N)^2$ . An integral of the form
\begin{equation*}
\iint (\Tr
  (W_1))^{\alpha_1}(\overline{\Tr(W_1)})^{\beta_1}\cdots (\Tr(W_p))^{\alpha_p}(\overline{\Tr
  (W_P)})^{\beta_p} \d U \d V
\end{equation*}
 is non-zero (modulo 
${\rm O}\big(\frac 1 N\big)$) if and only if it can be written in the form 
\begin{equation*}\iint _{\cU(N)^2} |\Tr (W_1)|^{2\alpha_1} \cdots |\Tr (W_p)|^{2\alpha_p} \d U
\d V\end{equation*}
 in which case it is equal to 
$\alpha_1!\cdots \alpha_p! j_1^{\alpha_1} \cdots j_p^{\alpha_p}$, with $j_i=
j(w_i)$ the number of cyclic rotations of the word $w_i$, that are leaving
$w_i$ invariant.

Consequently, $\Tr (W_1), \ldots, \Tr(W_s)$ have the asymptotic moment
distribution\break (as $N\to \infty$) of $\sqrt{j_1} Z_1, \ldots, \sqrt{j_s} Z_s$,
where $Z_1, \ldots, Z_s$ are independent normal gaussian variables. 
\end{thm}


\begin{proof}
We rewrite the formula from the preceding section as follows.
\begin{equation}\label{sase}
\begin{split}
& \iint_{\cU(N)^2} \Tr(W_1) \cdots \Tr(W_p) \d U \d V  \\
& =
\sum_{\beta_1\cdots\beta_{n_1}\beta_{n_1+1}\cdots \beta_{n_s}=1}
\iint a_{\beta_1\beta_2}^{(1)}\cdots 
a_{\beta_{n_1}\beta_1}^{(n_1)}\cdots 
a_{\beta_{n_{p-1}+1}\beta_{n_{p-1}+2}}^{(n_{p-1}+1)}
\cdots \\
& \qquad \cdots
a_{\beta_{n_p}}\beta_{n_{p-1}+1}^{(n_p)} \d U \d V 
\end{split}
\end{equation}
where the symbols $a^{(1)} \cdots a^{(n_s)}$ belong to the set
$\{U,V,U^*,V^*\}$. Here $n_i-n_{i-1}$ is the lenght of the word
$w_i$, $i=1,2,\ldots,s$.

Denote by $\tilde{U}$ the set of all symbols $a^{(i)}$ that
are equal to the letter $u$, and similarly for $\tilde{U^*},
\tilde{V},\tilde{V^*}$. Because of \cite{Petz}, unless 
$\card \tilde{U} = \card \tilde{U}^* = n$,
$\card \tilde{V} = \card \tilde{V}^* = m$,
the integral is ${\rm O}(\frac{1}{N})$.

According to the formula in the preceding paragraph the integral
will be the summation over all bijection
$\sigma_u : \tilde{U} \to \tilde{U}^*$, 
$\theta_u : \tilde{U}^* \to \tilde{U}$, 
$\sigma_v : \tilde{V} \to \tilde{V}^*$, 
$\theta_v : \tilde{V}^* \to \tilde{V}$, 
of 
\begin{equation}
\label{cinci}
\sum_{\sigma_u ,\theta_u, \sigma_v, \theta_v}
W_{\sigma_u \circ \theta_v^{-1}}\cdot
W_{\sigma_v \circ \theta_v^{-1}}
N^{\sharp R(\sigma_u ,\theta_u, \sigma_v, \theta_v)}
\end{equation}
where $R(\sigma_u ,\theta_u, \sigma_v, \theta_v)$ is the
equivalence relation generated by 
$$
\sigma_u(i) + 1 \sim i, \quad \theta_u (i) \sim i +1, 
\quad \sigma_v(i) + 1 \sim i, \quad \theta_v(i)\sim i+1
$$
where by the operation $i+1$, we mean successively 
(when $i=n_1,n_2,\ldots,n_s$)
$$
n_1+1=1, \quad n_2 + 1 = n_1 + 1, \quad \ldots, \quad
n_s+1 = n_{s-1} + 1
$$
(coresponding to the cycle $(1,\ldots, n_1)(n_1+1,\ldots,n_2)(n_{s-1}+1,
\ldots,n_s))$.

The equivalence relation has no singletons and hence
$N^{\sharp R(\sigma_u ,\theta_u, \sigma_v, \theta_v)}$
is at most $N^{n+m}$, while the term 
$W_{\sigma_u \circ \theta_v^{-1}}\cdot
W_{\sigma_v \circ \theta_v^{-1}}$
is of the order $N^k$, where $k\leq n+m$, with equality
if and only if $\sigma_u = \theta_v^{-1}$,
$\sigma_v = \theta_v^{-1}$.

This means that the only non zero terms will come from equivalence
relations of the form
$$
\sigma_u(i) + 1 \sim i, \quad \sigma_u^{-1} (i) \sim i +1, 
\quad \sigma_v(i) + 1 \sim i, \quad \sigma_v^{-1}(i)\sim i+1.
$$
To obtain the number of terms $m$ the sum (6) in this case
we need to determine the permutations $\sigma_u, \sigma_v$ for which
this equivalence relation has all the classes of two elements.

Denote by $\gamma$ the permutation with the property that
$\gamma|\tilde{U} = \sigma_u$, 
$\gamma|\tilde{U}^* = \sigma_u^{-1}$, 
$\gamma|\tilde{V} = \sigma_v$, 
$\gamma|\tilde{V}^* = \sigma_v^{-1}$.
Then $\gamma$ is an involution and the equivalence relation
is described as 
$$
\gamma(i) + 1 \sim i, \quad \gamma(i) \sim i+1.
$$
But $i \sim \gamma(i)+1$ and $i = (i-1)+1 \sim \gamma(i-1)$
and hence $\gamma(i-1)=\gamma(i)+1$ for all $i$ (since the
equivalence relations have only singletons). 
This means that if $i$ runs over the elements in a word 
$w_1$, then $\gamma(i)$ must run in the opposite direction
over the elements of the conjugate word $w_1^{-1}$.

In consequence, the integral in the statement is non zero
(modulo ${\rm O}(\frac{1}{N})$) only if it is of the form
$\iint_{\cU(N)} |\Tr(W_1)|^{2\alpha_1} \cdots
|\Tr(W_p)|^{2\alpha_p}\d u \d v$
and in this case the integral is equal to the number of possible
pairing between a word and cyclic rotation of its inverse.
This completes the proof.
\end{proof}

\section{A combinatorial formula for 
$\iint_{\cU(N)^2} |\Tr(W_1)|^2\cdots |\Tr(W_p)|^2\d U \d V$}

In this section we establish a formula that is specifically adapted
for integrals of products of absolute values of traces. Indeed
a positive answer for the Connes's embedding conjecture would
require the joint distribution of the variables
$|\Tr(W_1)|\cdots |\Tr(W_p)|$ as functions on $\cU(N)$.

Thus let $w_1,\ldots,w_p$ be reduced words in $F_2$, and let $X$
be the total set of symbols of 
$U,U^*,V,V^*$ in $|\Tr(w_1)|^2\cdots |\Tr(w_p)|^2$, where
each element in $X$ corresponds to a specific occurence of the
corresponding sysmbol in $|\Tr(w_1)|^2\cdots |\Tr(w_p)|^2$.

Assume there are $n$ occurences for the symbol $U$, $m$ occurences
for the symbol $V$, and hence that the set $X$ has $2(n+m)$ elements.
As in the preceding section $X$ is partitioned as 
$\tilde{U} \cup \tilde{U}^* \cup \tilde{V} \cup \tilde{V}^*$,
where $\tilde{U}, \tilde{U}^*, \tilde{V}, \tilde{V}^*$ are the set
of symbols of $u,u^*,v,v^*$ respectively in $X$.

Let $\Psi$ be the map which associates to each symbol $a$ in $X$,
which comes from a word $w_1$ its corresponding symbol $a^*$ in 
$w_1^{-1}$ and viceversa for a symbol $a$ in $w_1^{-1}$ it associates
the corresponding symbol $a^*$ in $w_1$.
Then $\Psi$ is an involution, $\Psi$ maps $\tilde{U}$ onto $\tilde{U}^*$
and $\tilde{V}$ onto $\tilde{V}^*$. Let $I$ be the map associating to
each symbol $a$ the successor of $\Psi(a)$ in the inverse word.

Let $S_{\tilde{U}}$ (and respectively $S_{\tilde{U}^*},
S_{\tilde{V}},S_{\tilde{V}^*}$) be the set of permutations of the
sets $\tilde{U}$ (respectively $\tilde{U}^*,\tilde{V},\tilde{V}^*$).
For each $\sigma_u \in S_{\tilde{U}}$, $\theta_u \in S_{\tilde{U}^*}$,
$\sigma_v \in S_{\tilde{V}}$, $\theta_v \in S_{\tilde{V}^*}$ let
$(\sigma_u, \theta_u, \sigma_v, \theta_v)$ be the concatanation of
these permutations to a permutation of $X$.

With the above notations we have:

\begin{prop}
Let $w_1,\ldots,w_p$ be words in $F_2$ and let $W_1,\ldots,W_p$
be the corresponding words as functions on $\cU(N)^2$.

Let $R$ be the following element in 
\begin{align*}
& \C(S_{\tilde{U}})\otimes
\C(S_{\tilde{U}^*})\otimes \C(S_{\tilde{V}})\otimes \C(S_{\tilde{V}^*}): 
\\
& R = \sum_{\scriptstyle \sigma_u \in S_{\tilde{U}}, \, \theta_u \in S_{\tilde{U}^*}
\atop \scriptstyle \sigma_v \in S_{\tilde{V}}, \, \theta_v \in S_{\tilde{V}^*}}
N^{\sharp (I\circ(\sigma_u ,\theta_u, \sigma_v, \theta_v))}
\sigma_u \otimes \theta_u \otimes \sigma_v \otimes\theta_v.
\end{align*}

Note that $R$ depends only of the cardinalities of the sets
$I(\tilde{U}) \cap U^*$,  
$I(\tilde{U}) \cap \tilde{V}$,  
$I(\tilde{U}) \cap \tilde{V}^*,\ldots,I(\tilde{V^*}) \cap \tilde{V}$.

Let $\Phi$ be the linear map from  
$\C(S_{\tilde{U}})\otimes
\C(S_{\tilde{U}^*})\otimes \C(S_{\tilde{V}})\otimes \C(S_{\tilde{V}^*})$
into $\C$ which associates to $\sigma_1\otimes \theta_1 \otimes
\sigma_2 \otimes \theta_2$ the number 
$W_{\sigma_1 \Psi^{-1}\theta_1 \Psi}^N
W_{\sigma_2 \Psi^{-1}\theta_2 \Psi}^N$.
Then
$$\iint_{\cU(N)^2} |\Tr(W_1)|^2\cdots |\Tr(W_p)|^2\d U \d V=\Phi(R).$$
\end{prop}  

\begin{proof}
Introduce an indexing  of the elements in $X$ so that the symbol
corresponding to a term $u_{\beta_i\beta_{i+1}}$ in a word $w_i$
to correspond to a $u_{\beta_{\overline{i+1}}\beta_{\overline{i}}^*}$ 
for a word in $w_1^{-1}$. Here we use the convention that the
elements in $w_i$ have indexing after $\beta_1,\beta_2,\ldots$.
Then compating the integral 
$$
\iint_{\cU(N)^2} |\Tr(W_1)|^2\cdots |\Tr(W_p)|^2\d U \d V
$$
will amount to compute integrals of the form
$$
\int_{\cU(N)} u_{\beta_s\beta_{s+1}}\cdots 
u_{\beta_{\overline{r}}\beta_{\overline{r-1}}}\cdots
u_{\beta_{\overline{s+1}}\beta_{\overline{s}}}^*\cdots
u_{\beta_{r-1}\beta_{r}}^* \d U. 
$$
Thus here 
$$U=\{\ldots,s,\ldots,\overline{r}\}, \quad
U^* = \{\ldots,\overline{s+1},\ldots,r-1,\ldots\}.$$
Then the map $\Psi$ will map $s,\overline{r}$ onto
$\overline{s+1}$, $r-1$ respectively and I will map
$s,r$ into $\overline{s},\overline{r}$ respectively.

We use the formula from Section 2 and we have the sum over
permutations $\sigma_u$ of the symbols $\{\ldots,s,\ldots,\overline{r}\}$,
and permutations $\theta_u$ of the symbols 
$\{\ldots,\overline{s+1},\ldots,r-1,\ldots\}$.
Hence the corresponding equivalence relation corresponding
to these permutations (and the similar permutations for $\theta$
will be exactly $$\overline{s} \sim \sigma_u(s),\quad
\overline{r}\sim \sigma_u(r)\quad \mbox{and}\quad \overline{s+1} \sim
\theta_u(s+1), \quad r-1 \sim \theta_u (\overline{r}-1),$$ which
is exactly the equivalence relation corresponding to
$I\circ (\sigma_u ,\theta_u, \sigma_v, \theta_v)$. This completes
the proof.
\end{proof}

\section{An example for the computation of 
$\iint_{\cU(N)} 
|\Tr(U^{p^1_1}V^{\epsilon_1} U^{p^1_2} V^{\epsilon^1_2} \cdots
U^{p^1_{s_1}} V^{\epsilon^1_{s_1}})|^2 \cdots
|\Tr (U^{p_{1}^t} V^{\epsilon_{1}^t}\cdots 
U^{p_{s_t}^t} V^{\epsilon_{s_t}^t})|^2 
\d U \d V$}

We apply the algorithm in the preceding section for the calculation
of a product  of words in which between powers of $u$ of degree
at least 3 are intercolated powers of $v$ of degree $\pm 1$.
We will describe the structure of the element $R$ in such a 
case since $\Psi$ is easy to be descried in this situation.

\begin{prop}
For the integral, $|p_a^s|\geq 3$, $\epsilon_a^s = \pm 1$,
$$\iint_{\cU(N)} 
|\Tr(U^{p^1_1}V^{\epsilon_1} U^{p^1_2} V^{\epsilon^1_2} \cdots
U^{p^1_{s_1}} V^{\epsilon^1_{s_1}})|^2 \cdots
|\Tr (U^{p_{1}^t} V^{\epsilon_{1}^t}\cdots 
U^{p_{s_t}^t} V^{\epsilon_{s_t}^t})|^2 
\d U \d V$$ 
the element $R$ is described as follows:

Let $n$ be the total number of $u$'s and $m$ the total number
of $v$'s.

The structure of the element $R$, viewed as an element
of $\C(S_n)\otimes \C(S_n)\otimes \C(S_m)\otimes \C(S_m)$
is described as follows:

Since $p_i^a \geq 1$ we have that $n\geq m$.
Then we have a set
$$
X = \{1_x,2_x,\ldots, (n-m)_x\} \cup
\{1_y,2_y,\ldots, m_y\} = X_0 \cup Y$$
and the first factor $S_n$ is identified with $S(x)$
(permutations of $X$), while the second $S_n$ is identified
with $S(\overline{X})$ where
$$
\overline{X} = \{\overline{1}_x,\overline{2}_x,\ldots, \overline{(n-m)}_x\} \cup
\{\overline{1}_y,\overline{2}_y,\ldots, \overline{m}_y\} = \overline{X}_0 \cup \overline{Y}.
$$
We consider also two sets
$A = \{1_a,\ldots,m_a\}$ and 
$\overline{A} = \{\overline{1}_a,\ldots,\overline{m}_a\}$. Then the first
factor $S_m$ is identified with $S(A)$, while the second with
$S(\overline{A})$.

The map $I$ acts on $X_0\cup \overline{X}_0$ by mapping $i_x$
into $\overline{i}_x$ and $\overline{i}_x$ into $i_x$, while on the set
$Y \cup \overline{Y}$, $I$ maps $i_y$ into $i_a$ and
$\overline{i}_y$ into $\overline{i}_a$ (or $i_y$ into $\overline{i}_a$ and
$\overline{i}_y$ into $i_a$). $I$ is an involution.
Then
$$
R = \sum_{\scriptstyle \sigma,\overline{\sigma} \in S(X),S(\overline{X}), 
\atop \scriptstyle \theta,\overline{\theta} \in  S(A),S(\overline{A})}
N^{\sharp (I\circ(\sigma,\overline{\sigma},\theta,\overline{\theta}))}
\sigma\otimes\overline{\sigma}\otimes\theta\otimes\overline{\theta}.
$$
\end{prop}

\begin{proof}
This follows by identifying the sets $\tilde{U}$, $\tilde{U}^*$
from the preceding proposition with the sets $X$, $\overline{X}$,
while $\tilde{V}$, $\tilde{V}^*$ are identified with 
$Y$, $\overline{Y}$.
\end{proof}

\begin{rem}
The map $\Psi$ can be explicitely describe as a map from $X$
onto $\overline{X}$, in terms of the alternating signs in 
$p_1^i,\epsilon_1^i,\ldots,p_s^i,\epsilon_s^i$, while on the set $Y$
it simply maps $i_a$ into $\overline{i_a}$.
\end{rem}

\section*{Acknowledgements}
We are indebted to Professors F.~Goodman and C.~Frohman
for several disscutions during the elaboration of this paper.
The author is grateful to the EPSRC and to the Deparment
of Mathematics of the Cardiff University, where part of this
work was done.

\section*{Note added while editing the manuscript}
We have been informed that R.~Speicher, P. Sniady  and J.~Mingo
have obtained independently in a joint paper
in preparation the same theorem as our result in Section 4.

\end{document}